\documentclass[11pt,a4paper]{amsart}

\swapnumbers \theoremstyle{plain}
\newtheorem{thm}{Theorem}[section]
\newtheorem{lem}[thm]{Lemma}

\newtheorem{lem-defn}[thm]{Lemma and definition}
\newtheorem{prop-defn}[thm]{Proposition and definition}
\newtheorem{cor-defn}[thm]{Corollary and definition}

\newtheorem{cor}[thm]{Corollary}

\theoremstyle{definition}

\theoremstyle{definition}

\newcommand{\T}{\mathbb{T}}

\newcommand{\C}{\mathbb{C}}
\newcommand{\D}{\mathbb{D}}

\DeclareMathOperator{\tr}{tr}
 \DeclareMathOperator{\im}{Im}
\DeclareMathOperator{\ke}{Ker} \DeclareMathOperator{\rea}{Re}
\DeclareMathOperator{\diam}{diam}

\numberwithin{equation}{section}

\usepackage{soul,cancel,url}
\usepackage{amsfonts}
\usepackage{hyperref}

\usepackage{color}
\usepackage{amsmath}
    \usepackage{amsxtra}
    \usepackage{amstext}
    \usepackage{amssymb}
   \usepackage{latexsym}
   \usepackage{graphicx}

\title {On holomorphic matrices on bordered Riemann surfaces}
\author{J\"urgen Leiterer}
\address{Institut f\"ur Mathematik \\
Humboldt-Universit\"at zu Berlin \\Rudower Chaussee 25\\D-12489 Berlin , Germany}
\email{leiterer@mathematik.hu-berlin.de}



\begin{document}

\begin{abstract} Let $\D$ be the unit disk. Kutzschebauch and Studer \cite{KS} recently proved that, for each continuous map $A:\overline D\to \mathrm{SL}(2,\C)$, which is holomorphic in $\D$, there exist continuous maps $E,F:\overline \D\to \mathfrak{sl}(2,\C)$, which are holomorphic in $\D$, such that  $A=e^Ee^F$. Also they asked if this extends to arbitrary compact bordered Riemann surfaces. We prove that this is possible.
\end{abstract}

\maketitle
\section{Introduction}

Let $\overline X$ be a  compact bordered Riemann surface\footnote{In the sense of \cite[II.3A]{AS}, which includes that $\overline X$ is connected. For example, $\overline X$ can be the closure of a bounded smooth domain $X$  in the complex plane $\C$.},  and let $X$ be the interior of $\overline X$.  Denote by $\mathrm{SL}(2,\C)$ the group of complex $2\times 2$ matrices with determinant $1$, and by $\mathfrak{sl}(2,\C)$  its Lie algebra of complex $2\times 2$ matrices with trace zero.
We prove the following.
\begin{thm}\label{Theorem 2} Let $A:\overline X\to \mathrm{SL}(2,\C)$ be a continuous map, which is holomorphic in $X$. Then there exist continuous maps $E,F:\overline X\to\mathfrak{sl}(2,\C)$, which are holomorphic in $X$, such that $A=e^{E}e^{F}$ on $\overline X$.
\end{thm} Let $\overline \D$ be the closed unit disk in $\C$.
For $\overline X=\overline \D$, Theorem \ref{Theorem 2} was recently proved by
Kutzschebauch and Studer \cite[Theorem 2]{KS}. In \cite{KS} also the question is asked if  Theorem \ref{Theorem 2} is true in general, and it is noted that there is some problem to adapt in a straightforward way the proof of \cite{KS}  to the general case. The problem is that $\overline X$ need not be simply connected. Our proof of Theorem \ref{Theorem 2} is nevertheless some adaption of that proof in \cite{KS}.

Let $\mathcal A(\overline X)$ be the algebra of complex-valued functions which are continuous on $\overline X$ and holomorphic in $X$. The first step in our proof of Theorem \ref{Theorem 2} is the following.

\begin{lem}\label{Bass2}Let $a,b\in \mathcal A(\overline X)$ with $\{a=0\}\cap \{b=0\}=\emptyset$ and, moreover, $\{a=0\}\not=\overline X$. Then there exists $h\in \mathcal A(\overline X)$ such that $b+ga=e^h$.
\end{lem}

This in particular implies that the Bass stable rank of $\mathcal A(\overline X)$ is one.\footnote{For $\overline X=\overline\D$, this is well-known \cite{JMW}. I do not know if this is already known for non-simply connected $\overline X$.}

That the  Bass stable rank of $\mathcal A(\overline \D)$ is one, is an important ingredient of the proof of Theorem \ref{Theorem 2}  given  in \cite{KS} for $\overline X=\overline \D$. As pointed out there, this makes it possible to limit to  matrices of the form $\big(\begin{smallmatrix}a&b\\c&d\end{smallmatrix}\big)$ with $\{a=0\}=\emptyset$. In the same way, Lemma \ref{Bass2} makes it possible
to limit to   matrices of the form $\big(\begin{smallmatrix}e^h&b\\c&d\end{smallmatrix}\big)$, and, for matrices of this form, it is possible to adapt the proof  from \cite{KS}  to the case of non-simply connected   $\overline X$.

Let $\mathrm{M}(2,\C)$  be the algebra of all complex $2\times 2$ matrices, and $\mathrm{GL}(2,\C)$ the group of its invertible elements. Then, in the same way as in \cite[Corollary 1]{KS},  the following corollary
 can be deduced from Theorem \ref{Theorem 2}.

\begin{cor}\label{26.8.20}  Let $A:\overline X\to \mathrm{GL}(2,\C)$ be continuous on $\overline X$, holomorphic in $X$, and null-homotopic. Then there exist continuous maps $E,F:\overline X\to \mathrm M(2,\C)$, which are holomorphic in $X$, such that  $A=e^{E}e^{F}$ on $\overline X$.
\end{cor}

 The study of the question ``how many exponentials factors are necessary to represent a  given holomorphic matrix'' was started by
Mortini and Rupp \cite{MR2}. In the case of an invertible $2\times 2$ matrix with entries from $\mathcal A(\overline \D)$,  they proved   that four exponentials are sufficient \cite[Theorem 7.1]{MR2}. Then Doubtsov and  Kutzschebauch \cite[Proposition 3]{DK} improved this to three exponentials. Eventually Kutzschebauch and Studer obtained that two exponentials are sufficient, which cannot be further improved, by an example Mortini and Rupp \cite[Example 6.4]{MR2}. This example   shows that,  under the hypotheses of Theorem \ref{Theorem 2} or Corollary \ref{26.8.20},  in general there does not exist a continuous $B:\overline X\to \mathrm{M}(2,\C)$ with $A=e^B$. As noted in
\cite{DK}, to find such $B$ with values in $\mathfrak{sl}(2,\C)$ is impossible already by the fact that not every matrix in $\mathrm{SL}(2,\C)$ has a logarithm in $\mathfrak{sl}(2,\C)$.

\section{A sufficient criterion for the existence of a logarithm}

A matrix $\Phi\in\mathrm M(2,\C)$ will be often considered as the linear operator in $\C^2$ defined  by multiplication from the left by $\Phi$ (considering the vectors in $\C^2$ as column vectors). The kernel and the image of this operator will be denoted  by $\ke \Phi$ and  $\im \Phi$, respectively.
For $\Phi\in \mathrm{M}(2,\C)$  and $\lambda\in \C$ we often write $\lambda-\Phi$ instead of $\lambda I-\Phi$. A matrix $\Phi\in\mathrm M(2,\C)$ will be called a {\bf projection}, if it is a linear projection as an operator, i.e., if $\Phi^2=\Phi$.

\begin{lem}\label{2.10.20--} Let $X$ be a topological space and let $B:X\to \mathrm{SL}(2,\C)$ be continuous.
Suppose there exists a continuous complex-valued function $\lambda$ on $X$ such that, for all $\zeta\in X$,
\begin{itemize}
\item[(a)] $e^{\lambda(\zeta)}$ is an eigenvalue of $B(\zeta)$,
\item[(b)] $e^{\lambda(\zeta)}\not=e^{-\lambda(\zeta)}$.
\end{itemize}
Then there exists a uniquely determined map $F:X\to\mathfrak{sl}(2,\C)$ such that $B=e^F$ on $X$ and, for all $\zeta\in X$, $\lambda(\zeta)$  is an eigenvalue of $F(\zeta)$. This map is  continuous. If $X$ is a complex space\footnote{By a complex space we mean a {\em reduced} complex space in the terminology of \cite{GR}, which is the same as an analytic space in the terminology of \cite{C} and \cite{L}. For example, each Riemann surface is a complex space.} and $B$, $\lambda$ are holomorphic, then $F$ is even holomorphic.
\end{lem}
\begin{proof} {\em Existence:} Since  $e^{\lambda(\zeta)}$ is an eigenvalue of $B(\zeta)$ and $\det B(\zeta)=1$, $e^{-\lambda(\zeta)}$ is the other eigenvalue of $B(\zeta)$, which is distinct from $e^{\lambda(\zeta)}$, by condition (b).
Therefore
\[\C^2=\ke\big(e^{\lambda(\zeta)}-B(\zeta)\big)\oplus\ke\big(e^{-\lambda(\zeta)}-B(\zeta)\big)\quad\text{for all}\quad\zeta\in X,\] where ``$\oplus$'' means ``direct sum'' (not necessarily orthogonal). Let $P:X\to \mathrm M(2,\C)$ be the map which assigns to each $\zeta\in X$  the linear projection from $\C^2$ onto $\ke\big(e^{\lambda(\zeta)}-B(\zeta)\big)$ along $\ke\big(e^{-\lambda(\zeta)}-B(\zeta)\big)$. Then \begin{equation}\label{4.10.20}
B=e^{\lambda}P+e^{-\lambda}(I-P),
\end{equation}which implies that
\begin{equation}\label{30.9.20''}
P=\frac{1}{e^{\lambda}-e^{-\lambda}}B-
\frac{e^{-\lambda}}{e^{\lambda}-e^{-\lambda}}I.
\end{equation} This shows that $P$ is continuous on $X$ and, if $X$ is a complex space and $B,\lambda $ are holomorphic, then $P$ is even holomorphic on $X$. Now
\begin{equation}\label{4.10.20'}
F:=\lambda P-\lambda(I-P)
\end{equation}
 has the desired properties.

{\em Uniqueness:}  Let $\zeta\in X$ and $\Theta\in \mathfrak{sl}(2,\C)$ such that $e^\Theta=B(\zeta)$, and $\lambda(\zeta)$ is an eigenvalue of $\Theta$.  Then $\Theta$ and $B(\zeta)$ commute. By \eqref{30.9.20''} also $\Theta$ and $P(\zeta)$ commute. Therefore  $\Theta=\alpha P(\zeta)+\beta\big(I-P(\zeta)\big)$ for some numbers $\alpha,\beta\in \C$, which  then are the eigenvalues of $\Theta$, i.e.,  either $\alpha=\lambda(\zeta)$  and $\beta=-\lambda(\zeta)$, or $\alpha=-\lambda(\zeta)$ and $\beta=\lambda(\zeta)$. $\alpha=-\lambda(\zeta)$ and $\beta=\lambda(\zeta)$ is not possible, since otherwise, by condition (b) and by \eqref{4.10.20}, we would have
\[
e^{\Theta}=e^{-\lambda(\zeta)}P+e^{\lambda(\zeta)}\big(I-P(\zeta)\big)\not
=e^{\lambda(\zeta)}P(\zeta)+e^{-\lambda(\zeta)}\big(I-P(\zeta)\big)=B(\zeta).
\] Therefore  $\alpha=\lambda(\zeta)$ and $\beta=-\lambda(\zeta)$. Hence, by \eqref{4.10.20'},
\[
\Theta=\lambda(\zeta)P(\zeta)-\lambda{\zeta}(I-P(\zeta))=F(\zeta).
\]
\end{proof}

\section{Proof of Lemma \ref{Bass2} and Theorem \ref{Theorem 2}}\label{2.9.20}

In this section,  $\overline X$ is a compact bordered Riemann surface, where we assume (as always possible\footnote{One can take for $\widetilde X$ a noncompact open neighborhood of $\overline X$ in the  double of $\overline X$ (for the definition of the double of $\overline X$, see, e.g., \cite[II. 3E]{AS}).}) that $X$ is a bounded smooth domain in some larger open Riemann surface $\widetilde X$, and $\overline X$ is the closure of $X$ in $\widetilde X$. The boundary of $\overline X$ will be denoted by $\partial X$. If we speak about an {\bf open  subset} $U$ {\bf of} $\overline X$,  then we  always mean that $U$ is a subset of $\overline X$ which is open in the topology of $\overline X$ (and in general not open in  $\widetilde X$). For  $K\subseteq\overline X$,  let $\overline K$ be the closure of $K$ (in $\overline X$ or in $\widetilde X$).

If $U$ is an open subset of $\overline X$, then we denote by $\mathcal A(U)$ the algebra of continuous complex valued functions on $U$ which are holomorphic in $U\cap X$.

To prove Theorem \ref{Theorem 2}, we begin with the observation that
\begin{equation}\label{27.9.20n}\Theta e^\Phi\Theta^{-1}=e^{\Theta\Phi\Theta^{-1}}\quad\text{ for all } \Theta,\Phi\in\mathrm{GL}(2,\C).
\end{equation} This shows that conjugation does not change the number of exponential factors needed to represent a given matrix. As  in \cite{KS}, we will use  this observation several times.

Next we recall  some known facts (Lemma \ref{26.9.20-}, its Corollary \ref{24.9.20} and Lemma \ref{15.9.20'}), for completeness with proofs.

\begin{lem}\label{26.9.20-} Let $\alpha$ be a continuous $(0,1)$-form  on $\overline X$ (i.e., a continuous section over $\overline X$ of the holomorphic cotangential bundle of $\widetilde X$) which is $\mathcal C^\infty $ in $X$. Then there exists a continuous function $u:\overline X\to \C$ which is $\mathcal C^\infty $ in $X$ such that $\overline\partial u=\alpha$ in $X$.
\end{lem}
\begin{proof}  As observed by Forstneric, Forn{\ae}ss and Wold in \cite[Section 2, formula (8)]{FFW} (together with corresponding references), to solve the $\overline\partial$-equation on Riemann surfaces one can use the following know fact: There exists a  1-form, $\omega$, defined and holomorphic on $(\widetilde X\times\widetilde X)\setminus\Delta$, where $\Delta $ is the diagonal in $\widetilde X\times \widetilde X$, such that, if $h:U\to \C$ is a  holomorphic coordinate on some open set $U\subseteq \widetilde X$, then,
on $(U\times U)\setminus \Delta$, $\omega$ is of the form
\begin{equation}\label{4.10.20-}
\omega(\zeta,\eta)=\Big(\frac{1}{h(\zeta)-h(\eta)}+\theta_h(\zeta,\eta)\Big)dh(\zeta),\quad (\zeta,\eta)\in (U\times U)\setminus \Delta,
\end{equation}where $\theta_h$ is a holomorphic function on $U\times U$. Since $\overline X$ is compact, and $\alpha $ is continuous on $\overline X$, then it is clear that the  function $u:\widetilde X\to \C$
defined by
\[
u(\eta)=\frac{1}{2\pi i}\int_{\zeta\in X}\omega(\zeta,\eta)\wedge\alpha(\zeta) ,\quad\eta\in \widetilde X,
\]is continuous on $\widetilde X$. To prove that, in $X$, $u$ is $\mathcal C^\infty$ and solves the equation $\overline \partial u=\alpha$, we consider a point $\xi\in X$ and take an open neighborhoods $V$ and $U$ of $\xi$ such that $\overline V\subseteq U$, $U\subseteq X$ and there exists a holomorphic coordinate $h:U\to \C$ of $\widetilde X$. Further choose  a $\mathcal C^\infty$-function $\chi:\widetilde X\to [0,1]$ such that  $\chi=1$ in a neighborhood $\overline V$. Then $u=u_1+u_2+u_3$, where
\begin{align*}
&u_1(\eta)=\frac{1}{2\pi i}\int_{\zeta\in V}\omega(\zeta,\eta)\wedge\alpha(\zeta),\\
&u_2(\eta)=\frac{1}{2\pi i}\int_{\zeta\in X\setminus V}\chi(\zeta)\omega(\zeta,\eta)\wedge\alpha(\zeta),\\
&u_3(\eta)=\frac{1}{2\pi i}\int_{\zeta\in X\setminus V}(1-\chi(\zeta))\omega(\zeta,\eta)\wedge\alpha(\zeta).
\end{align*}Then $u_2$ and $u_3$ are holomorphic in $V$. Therefore it remains to prove that $u_1$ is $\mathcal C^\infty$  and  $\overline\partial u_1=\alpha$, in $V$. By \eqref{4.10.20-}, $u_1=u_1'+u_1''$, where
\[
u'_1(\eta)=\frac{1}{2\pi i}\int_{\zeta\in V}\frac{ dh(\zeta)\wedge\alpha(\zeta)}{h(\zeta)-h(\eta)}\quad\text{and}\quad u''_1(\eta)=\int_{\zeta\in V}\theta_h(\zeta,\eta)dh(\zeta)\wedge \alpha(\zeta).\] Since $\theta_h$ is holomorphic, $u''_1$ is holomorphic. Further
\[
(u'_1\circ h^{-1})(w)=\frac{1}{2\pi i}\int_{z\in h(V)}\frac{dz\wedge\big((h^{-1})^*\alpha\big)(z)}{w-z}\quad\text{for }w\in h(V).
\]Therefore, as well-known (see, e.g, \cite[Theorem 1.2.2]{H}), $u'_1\circ h^{-1}$ is $\mathcal C^\infty$ and $\overline \partial (u'_1\circ h^{-1})=(h^{-1})^*\alpha$, in $h(V)$, which implies that $u'_1$ is $\mathcal C^\infty$ and $\overline\partial u'_1=\alpha$, in $V$.
\end{proof}

\begin{cor}\label{24.9.20} Let $U_1, U_2$  be nonempty open subsets of $\overline X$ with $U_1\cup U_2=\overline X$, and let  $f\in \mathcal A(U_1\cap U_2)$. Then there exist $f_1\in \mathcal A(U_1)$ and $f_2\in \mathcal A(U_2)$ with $f=f_1-f_2$ on $U_1\cap U_2$.
\end{cor}
\begin{proof}
For
 $K\subseteq \overline X$,  we denote by $\partial_{\overline X}K$ the boundary  of $K$ with respect to the topology of $\overline X$ (which is, in general, smaller than the boundary in $\widetilde X$). Since $U_1$ and $U_2$  are open subsets of $\overline X$ and $U_1\cup U_2=\overline X$, we have
\[\overline{U_1\setminus U_2}\cap \overline{U_2\setminus U_1} =\emptyset.\]Therefore  we can find
a $\mathcal C^\infty$ function $\chi:\widetilde X\to [0,1]$ with $\chi=1$ in an $\widetilde X$-neighborhood of $\overline{U_1\setminus U_2}$, and $\chi=0$ in an $\widetilde X$-neighborhood of  $\overline{U_2\setminus U_1}$.
Then we have well-defined continuous functions $c_1:U_1\to \C$ and $c_2:U_2\to \C$ which are $\mathcal C^\infty$ in $X\cap U_1$ and $X\cap U_2$, respectively,  such that
\[c_1=\begin{cases}(1-\chi) f&\text{on }U_1\cap U_2,\\
0&\text{on }U_1\setminus U_2,\end{cases}\quad \text{and}\quad c_2=\begin{cases}-\chi f&\text{on } U_1\cap U_2,\\
0&\text{on }U_2\setminus U_1.\end{cases}
\] Then
\begin{align}&\label{26.9.20+}f=c_1-c_2\quad\text{on}\quad U_1\cap U_2,\\
&\label{26.9.20++}\overline\partial c_1=-\overline\partial\chi f=\overline\partial c_2\quad\text{on}\quad X\cap U_1\cap U_2.
\end{align}Relation \eqref{26.9.20++}  shows that there is a well-defined continuous $(0,1)$-form on $\overline X$, $\alpha$, which is $\mathcal C^\infty$ in $X$,   such that
\begin{equation}\label{26.9.20+++}
\alpha =\overline\partial c_j\text{ on } X\cap U_j,\quad\text{for}\quad j=1,2.
\end{equation}By the preceding lemma, we can find a continuous function $u:\overline X\to \C$ which is $\mathcal C^\infty$ in $X$ such that $\overline\partial u=\alpha$ in $X$.
Set $f_j=c_j-u$, $j=1,2$. Then,  by \eqref{26.9.20+++}, $f_j\in\mathcal A(U_j)$ and, by \eqref{26.9.20+}, $f=f_1-f_2$ on $U_1\cap U_2$.
\end{proof}

\begin{lem}\label{15.9.20'} For each $a\in \mathcal A(\overline X)$, either $\{a=0\}=\overline X$ or $\partial X\cap \{a=0\}$ is nowhere dense in $\partial X$.
\end{lem}
\begin{proof}Assume $\partial X\cap \{a=0\}$ is not nowhere dense in $\partial X$. Then there exist $\xi\in \partial X$ and an open subset $U$ of $\overline X$ with $\xi\in U$ and $a\equiv 0$ on $U\cap \partial X$. Then (by defintion of a bordered Riemann surface), we have an open  subset $V$ of $\overline X$ with $\xi\in V$,  and a homeomorphism $\varphi:V\to \big\{z\in \C\,\big\vert\,\vert z\vert< 1,\,\im z\ge 0\big\}$,  which is biholomorphic from $V\setminus \partial X$ onto $\big\{z\in \C\,\big\vert\, \vert z\vert< 1,\, \im z>0\big\}$ and such that $\varphi(V\cap \partial X)=]-1,1[$. Then the continuous function $a\circ\varphi^{-1}$ is holomorphic in $\big\{z\in\C\,\big\vert\, \vert z\vert<1,\,\im z>0\big\}$ and has the real value $0$ on $]-1,1[$. Therefore, by the Schwarz reflection principle, there is a holomorhic function $\widetilde a$ on $\big\{z\in\C\,\vert\;\vert z\vert<1\big\}$ with
\begin{equation}\label{15.9.20}\widetilde a=a\circ\varphi^{-1}\quad\text{on}\quad\big\{z\in\C\,\big\vert\, \vert z\vert<1,\, \im z\ge 0\big\}.
\end{equation} Since $a=0$ on $\varphi^{-1}\big(]-1,1[\big)=V\cap \partial X$, from \eqref{15.9.20} we get $\widetilde a=0$ on $]-1,1[$. Therefore $\widetilde a=0$ on $\big\{z\in \C\,\big\vert\, \vert z\vert<1\big\}$. Again by \eqref{15.9.20} this implies that $a=0$ on $V\setminus \partial X$. Hence ($\overline X$ is connected)  $\{a=0\}=\overline X$.
\end{proof}

The first step in the proof of Lemma \ref{Bass2} is the following lemma.

\begin{lem}\label{2.9.20--} Let $a,b\in \mathcal A(\overline X)$ such that $\{ a=0\}\cap \{b=0\}=\emptyset$.
Then there exist finitely many closed subsets  $K_1,\ldots,K_\ell$ of $\overline X$ such that
\begin{align}
&\label{20.9.20n} K_j\cap K_k=\emptyset\quad\text{for all}\quad1\le j,k\le \ell\text{ with }j\not=k,\\
&\label{20.9.20n'}\{a=0\}\subseteq K_1\cup\ldots\cup K_\ell,
\end{align}
and, for some open disks $\D_1,\ldots,\D_\ell$ contained in $\C\setminus\{0\}$,
\begin{equation}\label{20.9.20n''} b(K_j)\subseteq \D_j\quad\text{for }j=1,\ldots,\ell.
\end{equation}
\end{lem}

\begin{proof} If $\{a=0\}=\emptyset$, the claim of the Lemma is trivial.
Therefore we may assume that $\{a=0\}\not=\emptyset$.

First let $\partial X\cap\{a=0\}=\emptyset$. Since $\overline X$ is compact and $\{a=0\}$ has no accumulation points in  $X$, and since $\{a=0\}\not=\emptyset$, then $\{a=0\}$ consists of a finite number of points $\xi_1,\ldots,\xi_\ell\in X$. Then $b(\xi_1)\not=0$, $\ldots$, $b(\xi_\ell)\not=0$, and $K_1:=\{\xi_1\}$, $\ldots$, $K_\ell:=\{\xi_\ell\}$  have the desired  properties.

Now let $\partial X\cap\{a=0\}\not=\emptyset$. Fix a metric $\rho(\cdot,\cdot)$ on $\widetilde X$. For a subset $K$ of $\widetilde X$ we denote by $\diam K$ the diameter of $K$ with respect to this metric. Since $\overline X$ is compact, $a,b$ are continuous and $\{a=0\}\cap\{b=0\}=\emptyset$, we have
\[\theta:=\min_{\zeta\in \overline X}\big(\vert a(\zeta)\vert+\vert b(\zeta)\vert\big)>0,\] and we can find $\varepsilon>0$ such that
\begin{equation}\label{11.9.20}
\big\vert b(\zeta)-b(\eta)\big\vert<\theta\quad\text{for all }\zeta,\eta\in \overline X\text{ with }\rho(\zeta,\eta)<\varepsilon.
\end{equation}

We call a set $\Lambda\subseteq \partial X$ a {\bf closed Interval} in $\partial X$ if there is a homeomorphic map $\psi$ from $[0,1]$ onto $\Lambda$.

Since $\overline X$ is  compact, $\partial X$ is the union of a finite number of pairwise disjoint Jordan curves.

{\bf Statement 1.} Let $\Gamma $ be one of these Jordan curves.  Then  there exists a finite number of closed intervals $\Lambda_1,\ldots,\Lambda_q$ in $\Gamma$  such that
\begin{align}
&\label{3.9.20''} \Lambda_j\cap \Lambda_k=\emptyset\quad\text{for }\quad 1\le j,k\le q\text{ with }j\not=k,\\
&\label{11.9.20-} \Gamma\cap \{a=0\}\subseteq \Lambda_1\cup\ldots\cup\Lambda_q,\\
&\label{20.9.20*} \Lambda_j\cap \{a=0\}\not=\emptyset\quad\text{for  }j=1,\ldots,q,\\
&\label{3.9.20'''} \diam(\Lambda_j)<\varepsilon\quad\text{for}\quad 1\le j\le q,
\end{align}
{\em Proof of Statement 1.} If $\Gamma\cap\{a=0\}=\emptyset$, the claim of the statement is trivial. Therefore we may assume that $\Gamma\cap\{a=0\}\not=\emptyset$.

Since $\Gamma$ is a Jordan curve, we have a homeomorphism $\phi$ from
$\T:=\{z\in \C\,\vert\, \vert z\vert=1\}$ onto $\Gamma$. Since $\{a=0\}\not=\overline X$, $\{a=0\}\cap \Gamma$ is nowhere dense in $\Gamma$ (Lemma \ref{15.9.20'}). Therefore we can find
 $0<t_1<t_2<\ldots<t_p<2\pi$ such that
 \begin{equation}\label{17.9.20}
 a\big(\phi(e^{it_\kappa})\big)\not=0\quad\text{for}\quad \kappa=1,\ldots,p,
 \end{equation} and
 \begin{multline}\label{17.9.20'}
\diam\phi\Big(e^{i[t_\kappa,t_{\kappa+1}]}\Big)<\varepsilon\text{ for }\kappa=1,\ldots,p-1,\text{ and}\\
\diam\Big(\phi\Big(e^{i[t_p,2\pi]}\Big)\cup\phi\Big(e^{i[0,t_1]}\Big)\Big)<\varepsilon.
 \end{multline}
By \eqref{17.9.20}, we can find $\sigma>0$ such that  $t_\kappa+\sigma <t_{\kappa+1}$ for $\kappa=1,\ldots,p-1$,  $t_p+\sigma<2\pi$, and
\begin{equation}\label{17.9.20''}
a\big(\phi(e^{it})\big)\not=0\text{ for }t_j\le t \le t_j+\sigma\text{  and }\kappa=1,\ldots,p.
\end{equation}
Define closed intervals in $\Gamma$, $\Delta_1,\ldots,\Delta_p$, by \begin{multline*}\Delta_\kappa=\phi\Big(e^{i[t_\kappa+\sigma,t_{\kappa+1}]}\Big)\text{ for }\kappa=1,\ldots,p-1, \text{ and}\\ \Delta_p=
\phi\Big(e^{i[t_p+\sigma,2\pi]}\Big)\cup\phi\Big(e^{i[0,t_1]}\Big).
\end{multline*}
Then it is clear that
\begin{equation}\label{16.9.20}\Delta_\kappa\cap \Delta_\lambda=\emptyset\quad\text{for all } \kappa,\lambda\in\{1,\ldots,p\}\text{ with }\kappa\not=\lambda,
\end{equation}  from \eqref{17.9.20'} it follows that
\begin{equation}\label{17.9.20-}
\diam \Delta_\kappa<\varepsilon\quad\text{for}\quad \kappa=1,\ldots,p,
\end{equation}and from \eqref{17.9.20''} it follows that
\begin{equation}\label{21.9.20}
\Gamma\cap\{a=0\}\subseteq \Delta_1\cup\ldots\cup \Delta_p.
\end{equation}

Let $\{\kappa_1,\ldots,\kappa_q\}$ be the set of all  $\kappa\in \{1,\ldots,p\}$ with $\Delta_\kappa\cap \{a=0\}\not=\emptyset$ (such $\kappa$ exist, as $\Gamma\cap\{a=0\}\not=\emptyset\}$), and define $\Lambda_j=\Delta_{\kappa_j}$ for $j=1,\ldots,q$.
Then \eqref{3.9.20''} is clear by \eqref{16.9.20}. \eqref{11.9.20-}  and \eqref{20.9.20*} hold by \eqref{21.9.20} and  the definition of the set $\{\kappa_1,\ldots,\kappa_q\}$. \eqref{3.9.20'''} is clear by \eqref{17.9.20-}. Statement 1 is proved.

From Statement 1 we obtain
a finite number of closed intervals $\Lambda_1,\ldots,\Lambda_r$ in $\partial X$  such that
\begin{align}
&\label{18.9.20-} \Lambda_j\cap \Lambda_k=\emptyset\quad\text{for }\quad 1\le j,k\le r\text{ with }j\not=k,\\
&\label{18.9.20--}\partial X\cap\{a=0\}\subseteq \Lambda_1\cup\ldots\cup\Lambda_r,\\
&\label{20.9.20**}\Lambda_j\cap\{a=0\}\not=\emptyset\quad\text{for }j=1,\ldots, r,\\
&\label{18.9.20----} \diam(\Lambda_j)<\varepsilon\quad\text{for}\quad j=1,\ldots, r.
\end{align}
By \eqref{18.9.20-} and \eqref{18.9.20----}, we can find  open subsets $U_j$ of $\overline X$, $j=1,\ldots,r$, with
\begin{align}
&\label{28.9.20-}\Lambda_j\subseteq U_j\quad\text{for}\quad 1\le j\le r,\\
&\label{18.9.20}\overline U_j\cap\overline U_k=\emptyset\quad \text{for all }1\le j,k\le r\text{ with }j\not=k,\\&\label{19.9.20}\diam\big(\overline U_j\big)<\varepsilon\quad\text{for }j=1,\ldots, r.
\end{align}Note that then, by \eqref{20.9.20**},
\begin{equation}
 \label{19.9.20'} U_j\cap \{a=0\}\not=\emptyset\quad\text{for }j=1,\ldots,r.
\end{equation} Set $K_j=\overline U_j$ for $j=1,\ldots, r$. Then, by \eqref{18.9.20--} and \eqref{28.9.20-},
\begin{equation}\label{28.9.20--}
\big\{a=0\big\}\cap \big(\partial X\cup K_1\cup \ldots\cup K_r\big)=
\big\{a=0\big\}\cap \big(K_1\cup \ldots\cup K_r\big).
\end{equation}

{\bf Statement 2.}
$N:=\{a=0\}\cap \big(\overline X\setminus \big(\partial X\cup K_1\cup\ldots\cup K_r\big)\big)$ is finite.

{\em Proof of Statement 2.} Assume $N$ is infinite. Since $\overline X$ is compact, then $N$ has an accumulation point $\xi\in \overline X$. Since $\{a=0\}$ is closed,  $\xi\in \{a=0\}$. As  $\{a=0\}\cap X$ is discrete in $X$, this implies that $\xi\in\partial X\cap \{a=0\}$ and further, by \eqref{18.9.20--} and \eqref{28.9.20-}, that $\xi\in  U_1\cup\ldots\cup U_r$. In particular, with respect to the topology of $\overline X$, $\xi$ is an inner point of $\partial X\cup K_1 \cup\ldots\cup K_r$, which is not possible, for $\xi$ is an accumultation point of $N$ and therefore, in particular, an accumulation point of $\overline X\setminus \big(\partial X\cup K_1\cup\ldots\cup K_r\big)$. Statement 2 is proved.

Let $\xi_{r+1},\ldots,\xi_\ell$  the distinct points of $N$, and define $K_{j}=\{\xi_{j}\}$ for $j=r+1,\ldots,\ell$. We claim that  $K_1,\ldots,K_\ell$ have the desired properties \eqref{20.9.20n}-\eqref{20.9.20n''}.

Indeed, \eqref{20.9.20n} follows from \eqref{18.9.20} and the fact that  $\xi_{r+1},\ldots,\xi_\ell$ are pairwise distinct and lie in $N$ and, hence, outside $K_1\cup\ldots\cup K_r$. By \eqref{28.9.20--},
\[
\{a=0\}\cap \big(\partial X\cup K_1\cup \ldots\cup K_r\big)\subseteq K_1\cup\ldots\cup K_r,
\] and, by definition of $K_{r+1},\ldots, K_\ell$,
\[
\{a=0\}\cap\Big(\overline X\setminus \big(\partial X\cup K_1\cup\ldots\cup K_r)\Big)=N=K_{r+1}\cup \ldots\cup K_\ell.
\] Together implies \eqref{20.9.20n'}. To prove \eqref{20.9.20n''}, we first note that by \eqref{19.9.20'} and the definition of $K_{r+1},\ldots,K_\ell$,
for each $j\in \{1,\ldots,\ell\}$, we have a point $\xi_j\in K_j$ with $a(\xi_j)=0$. Since, by definition of $\theta$,
$\vert b(\xi_j)\vert\ge \theta$, setting $\D_j=\Big\{z\in\C\,\big\vert\,\vert z-b(\xi_j)\vert <\theta\Big\}$ we obtain open disks $\D_1,\ldots,\D_\ell\subseteq \C\setminus\{0\}$. Since  $\diam K_j<\varepsilon$ for $j=1\ldots,\ell$ (for $1\le j\le r$ this holds by
\eqref{19.9.20}, and for $r+1\le j\le \ell$, we have $\diam K_j=0$), now \eqref{20.9.20n''} follows from  \eqref{11.9.20}. \end{proof}

{\bf Proof of Lemma \ref{Bass2}.}   If  $\{a=0\}=\emptyset$, we take a constant $C>0$ so large that $\rea(1+b(\zeta)/C)>0$ for all $\zeta\in X$,  and set $g=Ca^{-1}$. Then
\[
b+ga=b+C=C\big(1+b/C\big)=e^{\log C+\log(1+b/C)},
\]where $\log$ is the main branch of the logarithm.

Now let $\{a=0\}\not=\emptyset$. By Lemma \ref{2.9.20--}, we can find
finitely many closed subsets  $K_1,\ldots,K_\ell$ of $\overline X$ and open disks $\D_1,\ldots,\D_\ell$ in $\C\setminus \{0\}$ satisfying \eqref{20.9.20n}-\eqref{20.9.20n''}. Choose open subsets $W_1,\ldots,W_\ell$ of $\overline X$ such that
\begin{align}&\label{29.9.20}K_j\subseteq W_j\quad\text{for}\quad 1\le j\le \ell,\\
&\label{22.9.20'} W_j\cap W_k=\emptyset\quad\text{for all}\quad1\le j,k\le \ell\text{ with }j\not=k,\\
&\label{22.9.20} b(W_j)\subseteq \D_j\quad\text{for }j=1,\ldots,\ell.
\end{align}
Since $D_j\subseteq \C\setminus\{0\}$, we can find holomorphic functions $\log_j:\D_j\to \C$ with $e^{\log_j z}=z$ for  $z\in \D_j$. Set $W=W_1\cup\ldots\cup W_\ell$ and $V=\overline X\setminus\{a=0\}$. Then, by \eqref{29.9.20} and \eqref{20.9.20n'}, $V\cup W=\overline X$, and, by  \eqref{22.9.20'} and \eqref{22.9.20}, we can define $f\in \mathcal A(W)$ setting
$f=\log_j\circ b$ on $W_j$. Then
\begin{equation}\label{3.9.20-}
b=e^f\quad\text{on }W.
\end{equation}
Since $a\not=0$ on $V$ and $f\in \mathcal A(W)$, we have $f/a\in A(V\cap W)$. Therefore, by Corollary \ref{24.9.20}, we can find $v\in\mathcal A(V)$ and $w\in \mathcal A(W)$ with $f/a=v-w$, i.e.,
\[
f+aw=av\quad\text{on}\quad V\cap W.
\]
Therefore, we have a function $h\in \mathcal A(\overline X)$ with
\begin{equation}
\label{29.9.20''} h=f+aw\quad\text{on}\quad W.
\end{equation}
 The series $\sum_{\mu=0}^\infty\frac{a^\mu w^{\mu}}{\mu!}\frac{bw}{\mu+1}$
converges uniformly on the compact subsets of $W$  to some  $s\in \mathcal A(W)$, and, by \eqref{29.9.20''} and \eqref{3.9.20-}, we have \[e^h-b=e^{f+aw}-b=b e^{aw}-b=b(e^{aw}-1)\quad\text{on}\quad W.\] Together this implies that, on $V\cap W=W\setminus\{a=0\}$,
\[
\frac{e^h-b}{a}= \frac{b}{a}\sum_{\mu=1}^\infty\frac{a^\mu w^\mu}{\mu!}
=\frac{b}{a}\sum_{\mu=0}^\infty\frac{a^{\mu+1} w^{\mu+1}}{(\mu+1)!}=
\sum_{\mu=0}^\infty\frac{a^\mu w^\mu}{\mu!}\frac{bw}{\mu+1}=s.
\] Therefore, we have a function $g\in \mathcal A(\overline X)$ with
$\displaystyle g=\frac{e^h-b}{a}$ on $V$ and $g=s$ on $W$. Then, on $V=\overline X\setminus \{a=0\}$, it is clear that
\[b+ga=b+\frac{e^h-b}{a}a=e^h.\]  Since $\overline X\setminus\{a=0\}$ is nowhere dense in $\overline X$, it follows by continuity
that $b+ga=e^h$ on all of $\overline X$. \qed

{\bf Proof of Theorem \ref{Theorem 2}.} For $f\in\mathcal A(\overline X)$, we denote by $\rea f$ and  $\vert f\vert$ the functions $\overline X\ni \zeta\to \rea f(\zeta)$, and $\overline X\ni \zeta\to \vert f(\zeta)|$, respectively. By $\mathcal A^{\mathrm{SL}(2,\C)}(\overline X)$ and $\mathcal A^{\mathfrak{sl}(2,\C)}(\overline X)$ we denote the sets of  continuous maps from $\overline X$ to $\mathrm{SL}(2,\C)$ and $\mathfrak{sl}(2,\C)$, respectively, which are holomorphic in $X$.

Now let $A\in\mathcal A^{\mathrm{SL}(2,\C)}(\overline X)$ be given.

If $A\equiv I$ or $A\equiv -I$, the claim of Theorem \ref{Theorem 2} is trivial. Therefore  it is sufficient to consider the following three cases:
\begin{itemize}
\item[(I)] $A$ is of the form $\big(\begin{smallmatrix}a&b\\c&d\end{smallmatrix}\big)$ with $\{c=0\}\not=\overline X$,
\item[(II)] $A$ is of the form $\big(\begin{smallmatrix}a&b\\0&d\end{smallmatrix}\big)$ with  $\{b=0\}\not=\overline X$,
\item[(III)] $A$ is of the form $\big(\begin{smallmatrix}a&0\\0&d\end{smallmatrix}\big)$ where neither $\{a=1\}=\{d=1\}=\overline X$ nor $\{a=-1\}=\{d=-1\}=\overline X$.
\end{itemize}
By observation \eqref{27.9.20n}, Case (II) can be reduced to Case (I), since
\[
\begin{pmatrix}0&1\\1&0\end{pmatrix}\begin{pmatrix}a&b\\0&d\end{pmatrix}\begin{pmatrix}0&1\\1&0\end{pmatrix}^{-1}
=\begin{pmatrix}d&0\\b&a\end{pmatrix}.
\]
Consider Case (III). Since $\det A\equiv 1$, then $a\not=0$ and $d=a^{-1}$ on $\overline X$. Moreover, then   $\{a-a^{-1}=0\}\not=\overline X$, for otherwise we would have $\{a^2=1\}=\overline X$, i.e., either $\{a=1\}=\{d=1\}=\overline X$ or $\{a=-1\}=\{d=-1\}=\overline X$. As
\[
\begin{pmatrix}1&0\\1&1\end{pmatrix}\begin{pmatrix}a&0\\0&a^{-1}\end{pmatrix}\begin{pmatrix}1&0\\1&1\end{pmatrix}^{-1}=
\begin{pmatrix}a&0\\a-a^{-1}&a^{-1}\end{pmatrix},
\] this shows, again by  \eqref{27.9.20n}, that also Case (III) can be reduced to Case (I).

 So, we may assume that $A=\big(\begin{smallmatrix}a&b\\c&d\end{smallmatrix}\big)$ where $\{c=0\}\not=\overline X$. Since also  $\{c=0\}\cap \{a=0\}=\emptyset$ (the values of $A$ are invertible), then we can apply Lemma \ref{Bass2}, which gives $g,h\in \mathcal A(\overline X)$ with $a+gc=e^h$ on $\overline X$. Then
 \[
 \begin{pmatrix}1&g\\0&1\end{pmatrix}A\begin{pmatrix}1&g\\0&1\end{pmatrix}^{-1}=
 \begin{pmatrix}e^h&\ast\\\ast&\ast\end{pmatrix}.
 \]Therefore, again by observation \eqref{27.9.20n}, finally we see that $
A=\big(\begin{smallmatrix}e^h&b\\c&d\end{smallmatrix}\big)$ with $ h,b,c,d\in \mathcal A(\overline X)$
can be assumed.

The remaining part of the proof is an adaption of the proof given in \cite{KS} for
$\overline X=\overline \D$.
Chose  $\delta>0$ so large that, on $\overline X$,
\begin{align}&\label{1.10.20'}\rea\big(e^\delta+e^{h-\delta}d\big)>0,\\
&\label{1.10.20'''}\Big\vert(1+e^{h-2\delta}d)^2-4e^{-2\delta}-1\Big\vert
<\frac{1}{2},
\end{align} and define
\[E=\begin{pmatrix}h-\delta &0\\0& \delta-h\end{pmatrix}\quad\text{and}\quad B=\begin{pmatrix}e^\delta_{} &e^{\delta-h}b\\ e^{h-\delta}c&e^{h-\delta}d\end{pmatrix}.\]
Then
\begin{equation}\label{2.10.20-} E\in\mathcal A^{\mathfrak{sl}(2,\C)}(\overline X),\quad B\in\mathcal A^{\mathrm{SL}(2,\C)}(\overline X),\quad\text{and}\quad A=e^EB\text{ on }\overline X.
\end{equation}

 It follows from  \eqref{1.10.20'''} that $\log\big((1+e^{h-2\delta}d)^2-4e^{-2\delta}\big)$ is well-defined, where, since $\vert\log z\vert<1$ if $\vert z-1\vert<1/2$,
\begin{equation}\label{2.10.20'}
\big\vert\log\big((1+e^{h-2\delta}d)^2-4e^{-2\delta}\big)\big\vert<1\quad\text{on}\quad\overline X.
\end{equation}Since
\begin{equation*}
\frac{(\tr B)^2}{4}-1=\frac{e^{2\delta}}{4}\Big(\big(1+e^{h-2\delta}d\big)^2
-4e^{-2\delta}\Big),
\end{equation*} this implies that also $\log\big(\frac{(\tr B)^2}{4}-1\big)$ is  well-defined, where
\begin{equation}\label{6.10.20}
\log\Big(\frac{(\tr B)^2}{4}-1\Big)=2\delta-\log 4+\log\Big(\big(1+e^{h-2\delta}d\big)^2
-4e^{-2\delta}\Big)\quad\text{on}\quad\overline X.
\end{equation}
Set
\[
\varphi=\exp\bigg(\frac{1}{2}\log\Big(\frac{(\tr B)^2}{4}-1\Big)\bigg)\quad\text{on}\quad \overline X.
\]Then, by \eqref{6.10.20},
\begin{equation*}
\varphi=\exp\bigg(\delta-\frac{\log 4}{2}\bigg)\exp\bigg(\frac{1}{2}\log\Big(\big(1+e^{h-2\delta}d\big)^2-4e^{-2\delta}\Big)\bigg).
\end{equation*} Since $\vert e^z-1\vert<1$ if $\vert z\vert<1/2$ and therefore, by  \eqref{2.10.20'},
\[
\bigg\vert\exp\bigg(\frac{1}{2}\log\Big(\big(1+e^{h-2\delta}d\big)^2-4e^{-2\delta}\Big)\bigg)-1\bigg\vert<1,
\]this shows that
\begin{equation}\label{2.10.20+}
\rea\varphi>0\quad\text{on}\quad \overline X.
\end{equation}

Since $\varphi^2=\frac{(\tr B)^2}{4}-1$, we see  that, for each $\zeta\in\overline X$,
\[
\theta_+(\zeta):=\frac{\tr B(\zeta)}{2}+ \varphi(\zeta)\quad\text{and}\quad \theta_-(\zeta):=\frac{\tr B(\zeta)}{2}- \varphi(\zeta)
\]are the eigenvalues of $B(\zeta)$, where $\theta_+(\zeta)\not=\theta_-(\zeta)$ (as
$\varphi(\zeta)\not=0$). Since $\det B(\zeta)=1$ and therefore $\theta_-(\zeta)=\theta_+(\zeta)^{-1}$, it follows that $\theta_+(\zeta)\not=\theta_+(\zeta)^{-1}$ for all $\zeta\in \overline X$. Since, by \eqref{1.10.20'}, also $\rea(\tr B)>0$, it follows from \eqref{2.10.20+} that $\rea \theta_+>0$ on $\overline X$. Therefore $\lambda=\log\theta_+$ is well-defined. So, we have found a function $\lambda\in \mathcal A(\overline X)$ with the property that, for all $\zeta\in \overline X$, $e^{\lambda(\zeta)}$ ($=\theta_+(\zeta)$) is an eigenvalue of $B(\zeta) $ and $\lambda(\zeta)\not=-\lambda(\zeta)$ (as $\theta_+(\zeta)\not=\theta_+(\zeta)^{-1}$).
This implies by Lemma \ref{2.10.20--} that there exists $F\in \mathcal A^{\mathfrak{sl}(2,\C)}(\overline X)$ with $B=e^F$. By \eqref{2.10.20-} this completes the proof of Theorem \ref{Theorem 2}.


\begin{thebibliography}{FR10}

\bibitem {AS} L. Ahlfors and L. Sario, {\em Riemann Surfaces}, Princeton, 1960.

\bibitem{C} H. Cartan, {\em Espaces fibr\'es analytiques}, in: Symposium International de Topologia Algebraica.
Universidad National Aut\'onomica de Mexico and UNESCO. 97-121 (1958).

\bibitem{DK} E. Doubtsov and F. Kutzschebauch, {\em Factorization by elementary matrices, null-homotopy and products of exponentials for invertible matrices over rings}, Anal. Math. Phys. {\bf 9}, 1005-1018 (2019)

\bibitem {FFW} J. E. Forn{\ae}ss, F. Forstneri\v{c} and E. F. Wold, {\em Holomorphic approximation: the legacy of Weierstrass, Runge, Oka-Weil, and Mergelyan}, in: D. Breaz et al. (ed.), Advancements in complex analysis. From theory to practice. Cham: Springer. 133-192 (2020).


\bibitem{GR} H. Grauert and R. Remmert, {\em Coherent analytic sheaves}, Springer, 1984.

\bibitem{H} L. H\"ormander, {\em An introduction to complex analysis in several variables}, 3rd edn., North-Holland, 1990.

\bibitem{JMW} P. W. Jones, D. Marshall, T. H.  Wolff, {\em  Stable rank of the disc algebra}, Proc. Amer. Math. Soc. {\bf 96}, no. 4, 603–604 (1986).

\bibitem{KS} F. Kutzschebauch and L. Studer, {\em Exponential factorizations of holomorphic maps}, Bull. London Math. Soc. {\bf 51}, 995-1004 (2019).

\bibitem{L} S. \L ojasiewicz, {\em Introduction to complex analytic geometry}, Birkh\"auser, 1991.

\bibitem{MR2} R. Mortini and R. Rupp, {\em Logarithms and exponentials in the matrix algebra $\mathcal M_2(A)$}, Computer Methods Funct. Theory {\bf 18},  53-87 (2018)

\end{thebibliography}
\end{document}